\documentclass[11pt,leqno]{article}
\usepackage{amsfonts,amssymb,amsmath,amsthm,latexsym}
\usepackage{color,graphicx}
\usepackage[titletoc]{appendix}

\usepackage{amsfonts} 
\usepackage{latexsym}
\usepackage{graphicx}
\usepackage[utf8]{inputenc}
\usepackage[T1]{fontenc}
\newtheorem{theorem}{\sc Theorem}[section]
\newtheorem{lemma}{\sc Lemma}[section]
\newtheorem{proposition}{\sc Proposition}[section]

\newtheorem{remark}{\sc Remark}[section]

\begin{document}
\title{Biorthogonal functions for complex exponentials and an application to the controllability of the Kawahara equation via
a moment approach} 
\author{ Ademir F. Pazoto \thanks{Institute of Mathematics, Federal University of Rio de Janeiro,
UFRJ,P.O. Box 68530, CEP 21945-970, Rio de Janeiro, RJ, Brazil . E-mail: {\tt ademir@im.ufrj.br}} \and  Miguel D. Soto Vieira
\thanks{Institute of Mathematics, Federal University of Rio de Janeiro, UFRJ, P.O. Box
68530, CEP 21945-970, Rio de Janeiro, RJ, Brazil. E-mail: {\tt dariosv@gmail.com}} }
\date{}
\maketitle
\begin{abstract}
The paper deals with the controllability properties of the Kawahara equation posed on a periodic domain. We show that
the equation is exactly controllable by means of a control depending only on time and acting on the system through a given
shape function in space. Firstly, the exact controllability property is established for the linearized system through a Fourier expansion of solutions and the analysis of a biorthogonal sequence to a family of complex exponential functions. Finally, the local controllability
of the full system is derived by combining the analysis of the linearized system, a fixed point argument and some Bourgain smoothing properties of the Kawahara equation on a periodic domain.

\end{abstract}



\section{Introduction}
The study of wave phenomena arising in dispersive media is of broad scientific interest and pertains
to a modern line of research which is important both scientifically and for potential applications.
Progress in the development of mathematical models has made it possible to understand such phenomena
in quite distinct fields and to solve problems that come to the fore. Within this context, the Korteweg-de Vries equation
(KdV) has been derived as a model for the unidirectional propagation of nonlinear,
dispersive waves in an impressive array of physical situations. In most cases when it is derived from more complex systems,
the KdV equation appears in the form
$$ u_t + u_x + \varepsilon uu_x + \delta u_{xxx} = 0,$$
where the small positive parameters $\varepsilon$ and $\delta$ are related to a small-amplitude and a long-wavelength assumption, respectively.
The unknown $u$ is a real valued functions of the variables $x$ and $t$ and subscripts indicate partial differentiation.

Another relevant dispersive wave model is the Kawahara equation \cite{kawahara}, also referred as fifth-order
KdV equation. The Kawahara equation occurs in the theory of magneto-acoustic waves in a plasma and in the theory of shallow water waves
with surface tension. In order to balance the nonlinear effect, Kawahara took into account the higher order effect of
dispersion and established the following equation to describe solitary-wave propagation in media:
\begin{equation}\label{kawahara}
u_t + \gamma u_x + \alpha u_{xxx} + \beta u_{xxxxx} + \rho uu_x = 0.
\end{equation}
The parameters $\gamma, \alpha, \beta, \rho \in \mathbb{R}$ with $\beta\neq 0$, and $\alpha$ and $\beta$ represent the effect of dispersion.

There is a vast literature devoted to the study of water waves ranging from coastal engineering preoccupations to a very theoretical mathematical analysis of
the equations. For instance, a large body of literature has been concerned with the questions of existence, uniqueness and continuous dependence of solutions
corresponding to initial data. However, there are many issues still open that deserve further attention. 
In this work, the goal is to advance the study of the initial-boundary value problems exploring the dynamics of dispersive equations by using
mathematical analysis from the controllability point of view. Due to the rapid development of new mathematical tools, since the late 1980s control theory
of nonlinear dispersive wave equations have attracted a lot of attention.
Particularly, control properties of the KdV equation have been intensively studied and
significant progresses have been made. For a quite complete revision on the subject, we recommend the
works \cite{cerpa1} and \cite{rosier-zhang1}. 
In contrast, there are relatively few works on the Kawahara equation for its control theory (see, for instance, \cite{capistrano1,capistrano2,zhang1,zhang2}).

Without loss of generality, we assume that the parameters given in \eqref{kawahara} are such that $\gamma=\alpha=\rho=1$ and $\beta=-1$. Thus, our attention is given to the following control system described by the Kawahara equation posed on a periodic domain:
\begin{equation}\label{a01-1}
\begin{cases}
u_t-u_{5x}+u_{3x}+u_x+uu_x=f(x)v(t),   &\text{in} \ \ (0,T)\times (0,2\pi),\\
\partial_x^{j}u(t,0)=\partial_x^{j} u (t,2\pi), &\text{in} \ \ \ (0,T),\\
u(0,x)=u_0(x),&\text{in}\ \ \ (0,2\pi),
\end{cases}
\end{equation}
for $j=0, 1, 2, 3, 4$. The goal is to drive the initial data $u_0$ to rest by using a control $v(t)$, depending only on time and acting
on the system through a given function in space $f(x)$. This type of control is often used and sometimes called
lumped or bilinear.

To be more precise, considerations will be given to the following exact controllability problem:
\vglue 0.1 cm
{\it Given $T > 0$, an initial state $u_0$, a final state $u_1$ and a profile $f$ in a certain Hilbert space,
find an appropriate control $v\in L^2(0, T)$, so that system \eqref{a01} admits a solution $u$ which satisfies $u(T,x)=u_1(x)$.}

{\it If one can always find a control input to guide the system described by \eqref{a01} from
any given initial state to zero, then the system is said to be {\bf exactly controllable.}}

Since our system is time reversible, this property is equivalent to the
{\bf null-controllability property} which asserts that any initial state in a certain Hilbert space can be driven to zero in
time $T$.

\vglue 0.1 cm

In order to make more precise the tools we employ to study this question, we introduce some notations:
Given any $v\in L^2(0,2\pi)$
and $k\in\mathbb{Z}$, we denote by $\widehat{v}_k$ the $k-$Fourier
coefficient of $v$,
\[
\widehat{v}_k=\frac{1}{2\pi}\int_0^{2\pi} v(x) e^{-ikx} \, {\rm
d}x. \]

Then, for any $s\in \mathbb{R}$, we define the Hilbert space
\begin{equation}\label{eq:hsp}
H^s_p(0,2\pi)=\left\{ v=\sum_{k\in\mathbb{Z}} \widehat{v}_k
e^{ikx} \in L^2(0,2\pi)\, \left| \,\sum_{k\in\mathbb{Z}}
|\widehat{v}_k|^2(1+k^2)^s <\infty\right.\right\}
\end{equation}
endowed with the inner product
\begin{equation}\label{eq.innp}
(v,w)_s=\sum_{k\in\mathbb{Z}} \widehat{v}_k \overline{\widehat{w}_k}
(1+k^2)^s.\end{equation}
We denote by $\|\cdot\|_s$ the norm corresponding to the inner product given by \eqref{eq.innp}. Then, we consider the following operator
associated to the space variable:
\begin{equation}\label{A}
\begin{cases}
\vspace{2mm}(D(A),A), \mbox{ where } D(A)=H^5_p(0,2\pi) \mbox{ and }\\
A: D(A) \subset L^2_p(0,2\pi) \rightarrow L^2_p(0,2\pi), \mbox{ such that } Au=\partial_x^5u-\partial_x^3u-\partial_x u.
\end{cases}
\end{equation}

Taking the considerations above into account, we first address the controllability problem for the linearized system.
More precisely,
\begin{equation}\label{a01}
\begin{cases}
u_t-u_{5x}+u_{3x}+u_x=f(x)v(t),   &\text{in} \ \ (0,T)\times (0,2\pi),\\
\partial_x^{j}u(t,0)=\partial_x^{j} u (t,2\pi), &\text{in} \ \ \ (0,T),\\
u(0,x)=u_0(x),&\text{in}\ \ \ (0,2\pi),
\end{cases}
\end{equation}
for $j=0, 1, 2, 3, 4$.

Controllability properties of linear systems have been studied for a long time with the aid of Fourier techniques.
Concerning system \eqref{a01}, we employ Fourier series expansion to reduce the null control problem to a equivalent
{\it moment problem}, whose solution is given in terms of an explicit biorthogonal sequence to a family of exponential $(e^{\lambda_m t})_{m\in \mathbb{Z}}$
in $L^2(0,T)$. Here, $\lambda_m$ are the eigenvalues of the differential operator $A$ defined in \eqref{A}. We recall that a family
of functions  $(\phi_m)_{m\in \mathbb{Z}}\subset L^2(0,T)$ with the property that
$$\displaystyle\int_0^T\phi_m(t)e^{\overline{\lambda_m} t}dt=\delta_{m n}, \quad \forall\, m, n\in \mathbb{Z},$$
where $\delta_{m n}$ is the Kronecker symbol, is a biorthogonal sequence to $(e^{\lambda_m t})_{m\in \mathbb{Z}}$. In order to obtain this sequence, we introduce a family
$\Psi_m(z)$ of entire functions of exponential type (see, for instance, \cite{young}), such that $\Psi_m(i\lambda_n)=\delta_{m n}$.
Then, by applying Paley–Wiener Theorem we obtain $\phi_m$ as the inverse Fourier transform of $\Psi_m$. Each $\Psi_m$ is obtained
from a Weierstrass product $P_m$ multiplied by an appropriate function $M_m$ with rapid decay on the real axis. Such a method was used
for the first time by Paley and Wiener \cite{paley-wiener} and, in the context of control problems, by Fattorini and Russell \cite{fattorini1,fattorini2}.

Once such family $(\phi_m)_{m\in \mathbb{Z}}$ is given, the control $v(t)$ for \eqref{a01}
is obtained by considering a linear combinations of functions $\phi_m$. Indeed, if we consider $u_0(x)=\sum_{m\in \mathbf{Z}}\widehat{u}^0_{m}e^{imx}$
and $f(x)=\sum_{m\in \mathbf{Z}}\widehat{f}_{m}e^{imx},\,\widehat{f}_{m}\neq 0$, the Fourier expansions of $u_0$ and $f$, respectively, the function
\begin{align}\label{v-intro}
    v(t)=\sum_{m\in\mathbb{Z}}\frac{\hat{u}_{m}^{0}}{\hat{f}_m}e^{T\lambda_m}\phi_m\left(T-t\right),\ \ \ \ t\in(0,T),
\end{align}
is a control for \eqref{a01} in time $T$, if the series converges in $L^2(0,T)$. The convergence depends on
some uniform boundedness, with respect to $m$, of the the family $(\phi_m)_{m\in \mathbb{Z}}$ in $L^2(0,T)$, which
are obtained by applying Plancherel Theorem. In addition, some assumptions
on $f$ and $u_0$ are necessary. More precisely, let $f\in L^2(0,2\pi)$ be, such that
\begin{align}\label{u0,f}
f(x)=\sum_{k\in \mathbf{Z}}\widehat{f}_{k}e^{ikx},\,\, \mbox{with}\,\,\widehat{f}_{k}\neq 0, \,\,\forall\,k\in \mathbf{Z}.
\end{align}
Assuming \eqref{u0,f}, for a given constant $\beta >0$ define the space
\begin{equation}\label{H}
\mathcal{H}=\left\{h\in L^2_p(0,2\pi):\displaystyle\sum_{k\in\mathbb{Z}}\left|\frac{\widehat{h}_k}{\widehat{f}_k}\right|^2e^{\beta k^6}<\infty \right\}.
\end{equation}
If $u_0\in \mathcal{H}$ and $\widehat{f}_k$ satisfies \eqref{u0,f},
the convergence of \eqref{v-intro} holds in $L^2_p(0,2\pi)$ and $v(t)$ is a control for \eqref{a01}. We remark that the choice of the space $\mathcal{H}$ defined in \eqref{H} is
related to the form of the eigenvalues of the operator $A$ defined in \eqref{A} and the growth of $\phi_m$ in $L^2(0,T)$. Indeed, the eigenvalues of the state operator corresponding
to \eqref{a01} are given by $\lambda_m=- im(m^4+m^2-1)$ and $||\phi_m||_{L^2(0,T)}$ increases exponentially with $m$, i. e., $||\phi_m||_{L^2(0,T)}\leq ce^{\nu m^6 t}$,
where $c$ and $\nu$ are positive constants. The choice of the initial data in $\mathcal{H}$ compensates the growth
of $\phi_m$ and ensure the converge of \eqref{v-intro} in $L^2(0,T)$. When considering models in which the corresponding state operator has eigenvalues
with negative real part, we can take $\beta=0$ in \eqref{H}.


The technique we describe above was employed in the study of several control problems, being the pioneering articles of Fattorini and Russell
\cite{fattorini1,fattorini2} one of the most relevant examples in the context of scalar parabolic equations. This method
is very efficient in the one-dimensional space setting and has also been successfully
applied in \cite{micu-bugariu,glass,micu-ortega-pazoto}. In particular, our analysis was inspired by the results obtained
in \cite{micu-bugariu,micu-de teresa,micu-ortega-pazoto} of the which we borrow some ideas.

In order to prove the local controllability property for the full system \eqref{a01-1} we apply a fixed point argument
and the controllability result obtained for the linear system \eqref{a01}. At that point, we remark that
the Bourgain smoothing properties of the Kawahara equation obtained in \cite{hirayama} play a key role in our proof.


Concerning the Kawahara equation posed on a periodic domain, the internal controllability and the stabilization problems were studied in \cite{zhang1,zhang2}.
Particularly, in \cite{zhang2}, the authors use the same approach as that developed in \cite{laurent-rosier-zhang} to obtain the global exact control and global
exponential stability for periodic solutions in $H^s$,
for $s\geq 0$. Bourgain spaces associated to the Kawahara equation, propagation
of compactness and propagation of regularity for the linear Kawahara equation are three key ingredients in their proofs.
More recently, in \cite{flores-smith}, the authors establish local exact control and local exponential stability
of periodic solutions of fifth order Korteweg-de Vries type equations in  $H^s$, for $s > 2$. A dissipative term
is incorporated into the control which, along with a propagation of regularity property, yields a
smoothing effect permitting the application of the contraction principle. It is important to emphasize that the results obtained
in all papers mentioned above \cite{flores-smith,zhang1,zhang2} do not give an answer to control problem addressed here. Moreover,
they have been proved employing a different approach with a control input supported in a given open set $\omega\subset (0,2\pi)$. To the best of our knowledge,
the study we develop for the Kawahara equation has not been addressed in the literature yet. Moreover, the available results
do not give an immediate answer to it.

The remainder of this paper is organized as follows: in Section 2, we first give an equivalent characterization
of the controllability problem in terms of the moment problem. The next steps are devoted to the construction of a biorthogonal sequence and
to prove the controllability of the system \eqref{a01}. The local controllability of the full system is established in Section 3  and,
finally, in Section 4, we present some comments and open problems.

\section{The Linear System}

In this section we study the controllability properties of the system \eqref{a01-1}. We start by showing the equivalence between
the controllability and the moment problems. In order to do this, a result concerning the existence of solutions of \eqref{a01}
is needed.

\subsection{The Moment Problem}

Let us first present a well-posedness  result for system \eqref{a01}.

\begin{theorem}
Given any $T >0$,  $F \in L^1(0,T;L^2(0,2\pi))$ and $u^0 \in L^2(0,2\pi)$, there exists a unique weak solution $u \in C([0,T ]; L^2(0,2\pi))$ of the problem
\begin{equation}\label{a1}
\begin{cases}
u_t-u_{5x}+u_{3x}+u_x=F(t,x),  &\text{in  } \ \ (0,T)\times (0,2\pi),\\
\partial_x^ju(t,0)=\partial_x^ju(t,2\pi), &\text{in} \ \ \ (0,T),\\
u(0,x)=u_0(x), &\text{in}\ \ \ (0,2\pi),
\end{cases}
\end{equation}
for $j=0, 1, 2, 3, 4$.
\end{theorem}
\begin{proof}
According to \cite{zhang1}, the operator $A$ defined in \eqref{A} generates a group of isometries in $L^2_p(0,2\pi)$. Hence, the result follows from the semigroup theory.
\end{proof}
 Having the well-posedness of \eqref{a01} in hands, we can give now the characterization of the controllability property in terms
of a moment problem. We refer to \cite{avdonin,komornik,zab} for a detailed discussion of the subject.
\begin{theorem}
Let $T>0$, $f\in L^2(0,2\pi)$ and $u_0\in L^2_p(0,2\pi)$, such that
$$u_0(x)=\sum_{n\in\mathbb{Z}}\widehat{u}_{n}^{0}e^{inx}\quad\mbox{ and }\quad f(x)=\sum_{n\in\mathbb{Z}}\widehat{f}_ne^{inx}.$$
Then, there exists a control $v\in L^2(0,T)$ such that the solution $u$ of  \eqref{a01} verifies $u(T,x)=0$ if, and only if, $v\in L^2(0,T)$ satisfies
\begin{equation}\label{equiv}
    \widehat{f}_n\int_{0}^{T}v(T-s)e^{\lambda_ns}ds=-\widehat{u}_{n}^{0}e^{T\lambda_n},
\end{equation}
where $\lambda_n=- in(n^4+n^2-1)$ are the eigenvalues of the operator $A$ defined in \eqref{A}.
\end{theorem}
\begin{proof}
We consider the ``adjoint'' system
\begin{equation}\label{adj1}
\begin{cases}
\varphi_t-\varphi_{5x}+\varphi_{3x}+\varphi_x=0,  &\text{in  } \ \ (0,T)\times (0,2\pi),\\
\partial_x^{j}\varphi(t,0)=\partial_{x}^{j}\varphi(t,2\pi),&\text{in} \ \ \ (0,T),\\
\varphi(T,x)=\varphi_T(x), &\text{in}\ \ \ (0,2\pi),
\end{cases}
\end{equation}
for $j=0, 1,2, 3, 4$. If we multiply the equation in \eqref{a01} by $\overline{\varphi}$ and integrate for parts in $(0,T)\times(0,2\pi)$, we deduce that $v\in L^2(0,T)$ is a control for \eqref{a01} if, and only if, it verifies
\begin{equation}\label{controle}
    \int_{0}^{T}v(t)\int_{0}^{2\pi}f(x)\overline{\varphi}(t,x)dxdt=-\int_{0}^{2\pi}u_0(x)\overline{\varphi}(0,x)dx,
\end{equation}
for any solution $\varphi$ of \eqref{adj1}. Since $(e^{-inx})_{n\in\mathbb{Z}}$ is a basis for $L_p^2(0,2\pi)$, it is sufficient to check \eqref{controle} for solutions of \eqref{adj1} of the form $\varphi(t,x)=e^{(t-T)\lambda_n}e^{-inx}$, $n\in\mathbb{Z}$. Thus, it is straightforward to deduce that \eqref{equiv} holds.
\end{proof}

\subsection{A Biorthogonal Sequence}

This section is devoted to construct a biorthogonal sequence $(\Phi_m)_{m\in\mathbb{Z}}$ mentioned in the previous sections. By using Paley-Wiener Theorem,
it is obtained as the inverse Fourier transform of a family $\Psi_m$ of entire functions of exponential type, such that
$\Psi_m(i\lambda_n)=\delta_{mn}$, where $\delta_{mn}$ is the Kronecker symbol. Each $\Phi_m$ is obtained  from a Weierstrass product $P_m$ multiplied
by an appropriate function $M_m$ with rapid decay on the real axis. Therefore, for any $m\in \mathbb{Z}^*$, we first introduce the function
\begin{align}\label{funcion 1}
    P_m(z)=\displaystyle\prod_{n\in\mathbb{Z}^*, n\neq m}\left(1+\frac{iz}{\lambda_n}\right)\left(\frac{\lambda_n}{\lambda_n-\lambda_m}\right),
\end{align}
where $\lambda_m$ are the eigenvalues of the operator $A$ defined in \eqref{A}. Since $\lambda_{-m}=\overline{\lambda_m}$, we prove the following result:
\begin{lemma}\label{prop1}
$P_m$ is an entire function of the exponential type, such that
$$P_m(i\lambda_n)=\delta_{mn},\quad m\in\mathbb{Z}^*,$$
where $\delta_{mn}$ is the Kronecker symbol.
\end{lemma}
\begin{proof}
We obtain the result by analyzing the following products:
\begin{align}\label{EmQm}
    E_m(z)=\displaystyle\prod_{n\in\mathbb{Z}^*, n\neq m}\left|1+\frac{iz}{\lambda_n}\right| \ \ \text{and}\ \ \ Q_m=\displaystyle\prod_{n\in\mathbb{Z}^*, n\neq m}\left|\frac{\lambda_n}{\lambda_n-\lambda_m}\right|.
\end{align}
First, observe that, for any $z\in \mathbb{C}$,
\begin{align*}
   E_m(z)&=\displaystyle\prod_{n\in\mathbb{Z}^+, n\neq m}\left|1+\frac{iz}{\lambda_n}\right|\displaystyle\prod_{n\in\mathbb{Z}^-, n\neq m}\left|1+\frac{iz}{\lambda_n}\right|
   =\displaystyle\prod_{n\in\mathbb{N}^*, n\neq m}\left|1+\frac{iz}{\lambda_n}\right|\left|1+\frac{iz}{\overline{\lambda_n}}\right|\\
   &=\text{exp}\left(\displaystyle\sum_{n=1}^\infty \ln\left|1-\frac{z^2}{|\lambda_n|^2}+2iz\mathcal{R}\displaystyle\left(\frac{1}{\lambda_n}\right)\right|\right)
   =\text{exp}\left(\displaystyle\sum_{n=1}^\infty \ln\left|1-\frac{z^2}{|\lambda_n|^2}\right|\right).
\end{align*}
Since
\begin{align*}
    \displaystyle\sum_{n=1}^\infty \ln\left|1-\frac{z^2}{|\lambda_n|^2}\right|&\leq \displaystyle\sum_{n=1}^\infty \ln\left(1+2\frac{|z|^2}{|\lambda_n|^2}\right)\leq \displaystyle\sum_{n=1}^\infty \ln\left(1+2\frac{|z|^2}{n^2}\right)\leq \int_{0}^{\infty}\ln\left(1+2\frac{|z|^2}{x^2}\right)dx,\\
    &=\sqrt{2}\pi|z|,
\end{align*}
we get
\begin{equation}\label{Em}
E_m(z)\leq \text{exp}(\sqrt{2}\pi|z|).
\end{equation}

For $Q_m$ have that:
\begin{align*}
    &Q_m=\displaystyle\prod_{n\in\mathbb{Z}^*, n\neq m}\left|\frac{\lambda_n}{\lambda_n-\lambda_m}\right|=\frac{1}{2}\displaystyle\prod_{n\in \mathbb{N}^* , n\neq m}\frac{|\lambda_n|^2}{|\lambda_n-\lambda_m||\lambda_n+\lambda_m|} \\
    &=\frac{1}{2}\underbrace{\displaystyle\prod_{n\in\mathbb{N}^*, n\neq m}\frac{|\lambda_n|^2}{|\lambda_{n-m}||\lambda_{n+m}|}}_{Q_{m}^{1}}\underbrace{\displaystyle\prod_{n\in\mathbb{N}^*, n\neq m}\frac{|\lambda_{n-m}||\lambda_{n+m}|}{|\lambda_n-\lambda_m||\lambda_n+\lambda_m|}}_{Q_{m}^{2}}.
\end{align*}
Then, the next steps are devoted to estimate $Q_m^1$ and $Q_m^2$.
\begin{align*}
    Q_{m}^{1}&=\displaystyle\prod_{n\in\mathbb{N}^*, n\neq m}\frac{|\lambda_n|^2}{|\lambda_{n-m}||\lambda_{n+m}|}\leq \frac{|\lambda_1|^2|\lambda_2|^2\cdots|\lambda_{m-1}|^2|\lambda_{m+1}|^2\cdots|\lambda_{2m-1}|^2|\lambda_{2m}|^2|\lambda_{2m+1}|^2\cdots}{|\lambda_{m-1}|\cdots|\lambda_1||\lambda_{m+1}|\cdots|\lambda_{2m-1}|\displaystyle\prod_{n=1}^{\infty}|\lambda_n|\displaystyle\prod_{n=2m+1}^{\infty}|\lambda_n|}\\
    &\leq \frac{|\lambda_{2m}|}{|\lambda_m|}=\frac{|32m^5+8m^3-2m|}{|m^5+m^3-m|}\leq C,
\end{align*}
where $C$ is a positive constant.

To evaluate $Q_{m}^{2}$, we proceed as follows:
\begin{align*}
    &Q_{m}^{2}=\displaystyle\prod_{n\in\mathbb{N}^*, n\neq m}\frac{|\lambda_{n-m}||\lambda_{n+m}|}{|\lambda_n-\lambda_m||\lambda_n+\lambda_m|}=\displaystyle\prod_{n\in\mathbb{N}^*, n\neq m}\left(1+\frac{|\lambda_{n-m}||\lambda_{n+m}|-|\lambda_n-\lambda_m||\lambda_n+\lambda_m|}{|\lambda_n-\lambda_m||\lambda_n+\lambda_m|}\right)\\
    &=\displaystyle\prod_{n\in\mathbb{N}^*, n\neq m}\left(1+\frac{|\lambda_{n-m}\lambda_{n+m}|-|(\lambda_n-\lambda_m)(\lambda_n+\lambda_m)|}{|\lambda_n-\lambda_m||\lambda_n+\lambda_m|}\right)\\
    &\leq \displaystyle\prod_{n\in\mathbb{N}^*, n\neq m}\left(1+\frac{|\lambda_{n-m}\lambda_{n+m}-(\lambda_n-\lambda_m)(\lambda_n+\lambda_m)|}{|\lambda_n-\lambda_m||\lambda_n+\lambda_m|}\right)\\
    &=\exp\left(\sum_{n=1,n\neq m}^{\infty}\ln\left(1+\frac{|\lambda_{n-m}\lambda_{n+m}-(\lambda_n-\lambda_m)(\lambda_n+\lambda_m)|}{|\lambda_n-\lambda_m||\lambda_n+\lambda_m|}\right)\right)\\
    &\leq\exp\left(\sum_{n=1,n\neq m}^{\infty}\left(\frac{|\lambda_{n-m}\lambda_{n+m}-(\lambda_n-\lambda_m)(\lambda_n+\lambda_m)|}{|\lambda_n-\lambda_m||\lambda_n+\lambda_m|}\right)\right)\\
    &\leq\exp\left(\sum_{n=1,n\neq m}^{\infty}\frac{5m^8f(\frac{n}{m})+4m^6g(\frac{n}{m})+13m^4h(\frac{n}{m})}{\alpha(m,n)}\right),
\end{align*}
 where
\begin{align*}
\vspace{2mm}&f(t)=t^6-t^4+t^2, \quad g(t)=t^4+t^2,\quad h(t)=t^2,\\
&\alpha(m,n)=(n^4+n^3m+n^2m^2+nm^3+m^4+n^2+nm+m^2-1)\times\\
     &\qquad\qquad|n^4-n^3m+n^2m^2-nm^3+m^4+n^2-nm+m^2-1|.
     \end{align*}

In the remaining part of the proof $C$ will denote a positive constant that may change from one estimate to another, but it is independent of $m$.

Observe that the function $f(t)$ satisfy
 \begin{equation}
 f(t)\leq  \left\{\begin{array}{l}
t^2,\,\,\mbox{if}\quad 0\leq t\leq 1\\
t^6,\,\,\mbox{if}\quad t\geq 1.\nonumber
\end{array}\right.
 \end{equation}
Then, if $n\leq m$,
 \begin{align*}
     \sum_{n=1}^{m-1}\frac{5m^8f(\frac{n}{m})}{\alpha(m,n)}
     \leq5m^8 \sum_{n=1}^{m-1}\frac{\frac{n^2}{m^2}}{n^4}\leq 5m^6\sum_{n=1}^{m-1}\frac{1}{n^2}\leq 5m^6\int_{1}^{m-1}\frac{1}{t^2}dt=5m^6\frac{m-2}{m-1}\leq 5m^6.
 \end{align*}
 If $n\geq m$,
 \begin{align}\label{f}
     \sum_{n=m+1}^{\infty}\frac{5m^8f(\frac{n}{m})}{\alpha(m,n)}\leq 5m^8\sum_{n=m+1}^{\infty}\frac{\frac{n^6}{m^6}}{n^4(n-m)^4}\leq 5m^2\sum_{n=m+1}^{\infty}\frac{n^2}{(n-m)^4}\leq 5m^2\sum_{k=1}^{\infty}\frac{(k+m)^2}{k^4}\leq Cm^4.
 \end{align}
In what concerns the function $g(t)$, have that
 \begin{equation}
 g(t)\leq\left\{\begin{array}{l}
(t+1)^2,\,\mbox{if}\quad 0\leq t\leq 1,\\
2t^6,\,\qquad\mbox{if}\quad t\geq 1.
\end{array}\right.\nonumber
 \end{equation}
 When $n\leq m$,
 \begin{align*}
     \sum_{n=1}^{m-1}\frac{2m^6g(\frac{n}{m})}{\alpha(m,n)}&\leq 2m^6\sum_{n=1}^{m-1}\frac{(\frac{n}{m}+1)^2}{n^4}\leq 2m^4\sum_{n=1}^{m-1}\frac{(n+m)^2}{n^4}\leq 2m^4\sum_{n=1}^{m-1}\frac{n^2+2nm+m^2}{n^4}\\
     &\leq 2m^4\int_{1}^{m-1}\left(\frac{1}{t^2}+\frac{2m}{t^3}+\frac{m^2}{t^4}\right)dt\leq Cm^6.
 \end{align*}
 If $n\geq m$, we proceed as in \eqref{f}. In this case, we use the fact that $g(t)\leq 4t^6$, for $t\geq 1$. Finally, to estimate the term involving the function $h$, we also proceed as before using the following estimate:
 \begin{equation}
 h(t)\leq\left\{\begin{array}{l}
t^2,\quad\mbox{if}\quad 0\leq t\leq 1,\\
t^6,\quad\mbox{if}\quad t\geq 1.
\end{array}\right.\nonumber
 \end{equation}
 Combining the estimates above, we deduce that
 $$Q_m=Q^1_m Q^2_m\leq \exp(Cm^6).$$
From \eqref{EmQm}, \eqref{Em} and the above estimate we conclude the proof.
\end{proof}
\begin{remark}
Lemma \ref{prop1} remains valid if we consider the following linear equation equation associated to \eqref{kawahara}: $u_t + \gamma u_x + \alpha u_{xxx} - \beta u_{xxxxx}= 0$.
In fact,  the differential operator associated to the space variable is given by $A_1:=\beta\partial_x^5u-\alpha\partial_x^3u-\gamma\partial_x u:H^5_p(0,2\pi)\rightarrow L^2(0,2\pi)$, whose eigenvalues are
\begin{equation*}
    \lambda_k=-ik(\beta k^4+\alpha k^2-\gamma), \quad k\in\,\mathbb{Z}.
\end{equation*}
Hence, it may occur that not all eigenvalues are different. If we count only the distinct eigenvalues, we get a sequence $\{\lambda_k\}_{k\in \mathbb{I}}$, where $\mathbb{I}\subset \mathbb{Z}$ have a property of $\lambda_{k_1}\neq\lambda_{k_2}$ for any $k_1,k_2\in \mathbb{I}$. Then, for all $k_1\in \mathbb{Z}$, we define
\begin{align*}
    I(k_1)=\{k\in \mathbb{Z}: k(\beta k^4+\alpha k^2-\gamma)=k_1(\beta k_{1}^{4}+\alpha k_{1}^{2}-\gamma)\}
\end{align*}\label{Ik}
and $|I(k_1)|=m(k_1)$, which has the following properties:
\begin{itemize}
    \item $m(k_1)\leq 5$. This is a consequence of the fact that the polynomial $p(x)=x(\beta x^4+\alpha x^2 - \gamma)$ has a maximum of $5$ distinct roots.
    \item  $\lambda_k\rightarrow \pm \infty$, as $k\rightarrow\pm\infty$. Then, there exists $k^*\in \mathbb{N}$, such that $m(k)=1$, for all $|k|\geq k^*$.
\end{itemize}
To prove Lemma \ref{prop1}, we have assumed that $I(k_1)$ is a unitary set. This is due to the fact that, in the original model, we have assumed that
$\alpha=\beta=\gamma=1$. If this is not the case, we can also prove the result by using the same approach. Indeed, following the notation introduced in the proof
of the lemma, we have that
\begin{align*}
Q_{m}^{1}&=\displaystyle\prod_{n\in\mathbb{N}^*, n\notin I(m)}\frac{|\lambda_n|^2}{|\lambda_{n-m}||\lambda_{n+m}|}\\
&=\prod_{n=1}^{m_1-1}\frac{|\lambda_n|^2}{|\lambda_{n-m}||\lambda_{n+m}|}\prod_{m_1+1}^{m_2-1}\frac{|\lambda_n|^2}{|\lambda_{n-m}||\lambda_{n+m}|}
\cdots \prod_{m_5+1}^\infty\frac{|\lambda_n|^2}{|\lambda_{n-m}||\lambda_{n+m}|}.
\end{align*}
Then, proceeding in a similar way, we can estimate each term of the product above. For $Q_{m}^{2}$, we use a similar argument.
\end{remark}

From Lemma \ref{prop1} we obtain the following estimate for $P_m$, defined in \eqref{funcion 1}:
\begin{align}
   |P_m(z)|\leq \exp(C\pi(|z|+m^6),\nonumber
\end{align}
where $C$ is a positive constant. Consequently, on the real axis, it follows that
\begin{align}\label{estPmRe}
|P_m(x)|\leq \exp(C_1(|x|+m^6),
\end{align}
for some $C_1>0$.

The next proposition guarantees the existence of a entire function (of exponential type) which plays an important
role in the construction of the biorthogonal sequence. It is an appropriate multiplier that compensates the growth of $P_m$ on the real axis.
In order to prove the proposition, the following technical lemma is needed.

\begin{lemma}\label{a}
If $x\geq m^6$, then
\begin{align}
    \sum_{j=m^6}^{[x]}\ln \left|\frac{j}{x}\right|=-\int_{m^6}^{x}\frac{B(u)-m^6+1}{u}du,
\end{align}
where $B(u)=\#\{n: n\leq u\}$.
\end{lemma}

\begin{proof} Firstly, we remark that the function $B$ has the following properties:
\begin{itemize}
    \item If $j\leq u < j+1$, we have $B(u)=j$.
    \item If $[x]\leq u\leq x$, then $B(u)=[x]$ and $B(u)\geq x-1$.
\end{itemize}
Hence, we have that
\begin{align*}
    -\int_{m^6}^{x}\frac{B(u)}{u}du&=-\sum_{j=m^6}^{[x]-1}\int_{j}^{j+1}\frac{B(u)}{u}du-\int_{[x]}^{x}\frac{B(u)}{u}du\\
    &=-\sum_{j=m^6}^{[x]-1}\int_{j}^{j+1}\frac{j}{u}du-\int_{[x]}^{x}\frac{[x]}{u}du\\
    &=\sum_{j=m^6}^{[x]-1}j\ln\left|\frac{j}{j+1}\right|+[x]\ln\left|\frac{[x]}{x}\right|=\ln\left|\prod_{j=m^6}^{[x]-1}\frac{(j)^j}{(j+1)^j}\frac{([x])^{[x]}}{(x)^{[x]}}\right|\\
    &=\ln\left|\frac{(m^6)^{m^6}}{(m^6+1)^{m^6}}\frac{(m^6+1)^{m^6+1}}{(m^6+2)^{m^6+1}}\cdots\frac{([x]-1)^{[x]-1}}{([x])^{[x]-1}}\frac{([x])^{[x]}}{(x)^{[x]}}\right|\\
    &=\ln\left|\frac{(m^6)^{m^6-1}}{([x])^{m^6-1}}\prod_{j=m^6}^{[x]-1}\frac{j}{x}\right|=-\int_{m^6}^{x}\frac{m^6-1}{u}du+\sum_{j=m^6}^{[x]}\ln\left|\frac{j}{x}\right|.
\end{align*}
\end{proof}

As remarked above, Lemma \ref{a} allows us to prove the following result, inspired in \cite{ingham}:

\begin{proposition}\label{prop2}
For each $m\geq1$, there exists a function $M_m:\mathbb{C}\rightarrow \mathbb{C}$ and positive constants $K_1, K_2>0$, such that:
\begin{itemize}
    \item $M_m$ is a function of the exponential type,
    \item $|M_m(x)|\leq \exp(K_1(m^6-|x|)), \forall\, x\,\in \mathbb{R}$,
    \item $|M_m(i\lambda_m)|\geq \exp(-K_2m^6)$,
\end{itemize}
where $\lambda_m=-im(m^4+m^2-1)$ are the eigenvalues of the operator $A$ defined in \eqref{A}.
\end{proposition}
\begin{proof} We follow the ideas introduced in \cite{ingham} and define a function $M_m:\mathbb{C}\rightarrow\mathbb{C}$ as follows:
\begin{align}\label{funcion 2} M_m(z)=\prod_{n=m^3}^{\infty}\frac{\sin(\frac{z}{n^2})}{\frac{z}{n^2}}.
\end{align}
Since $\displaystyle\sum_{n=1}^\infty \frac{1}{n^2}<\infty$, the first property is a consequence of the following estimate:
\begin{align*}
 \prod_{n=m^3}^{N}\left|\frac{\sin(\frac{z}{n^2})}{\frac{z}{n^2}}\right| \leq \prod_{n=m^3}^{N}\exp\left({\left|\frac{z}{n^2}\right|}\right)= \exp({|z|\sum_{n=m^3}^{N}\frac{1}{n^2}})\leq \exp(C|z|),
\end{align*}
for some $C>0$.

To prove the second property, we proceed in two steps, as follows:

\vglue 0.2 cm

\noindent $\bullet$ If $|x|\leq m^6$, then
\begin{align*}
    |M_m(x)|=\prod_{n=m^3}^{\infty}\left|\frac{\sin(\frac{x}{n^2})}{\frac{x}{n^2}}\right|\leq 1\leq \exp(m^6-|x|).
\end{align*}
\noindent $\bullet$ If $|x|>m^6$, we apply Lemma \ref{a} to deduce that
\begin{align*}
    |M_m(x)|&=\prod_{n=m^3}^{\infty}\left|\frac{\sin(\frac{x}{n^2})}{\frac{x}{n^2}}\right|\leq \prod_{n=m^3}^{[|x|^{\frac{1}{2}}]}\frac{n^2}{|x|}=\exp\left(\sum_{n=m^3}^{[|x|^{\frac{1}{2}}]}\ln\frac{n^2}{|x|}\right)\leq\exp\left(\sum_{n=m^6}^{[|x|]}\ln\frac{n^2}{|x|}\right)\\
    &=\exp\left(-\int_{m^6}^{|x|}\frac{B(u)-m^6+1}{u}du\right).
\end{align*}
Since  $m^6\leq [|x|]$, from the estimate above, we obtain a positive constant satisfying
\begin{align*}
    |M_m(x)|&\leq\exp\left(-\int_{[|x|]}^{|x|}\frac{B(u)-m^6+1}{u}du\right)\leq \exp\left(-\int_{[|x|]}^{|x|}\frac{|x|-1-m^6+1}{u}du\right)\\
    &=\exp\left((m^6-|x|)\ln\frac{|x|}{[|x|]}\right)  \leq C\exp(m^6-|x|),
\end{align*}
where $C$ is a positive constant.

In what concerns the third property, we observe that $m^6\geq |\lambda_m|$, i. e., $\left|\displaystyle\frac{\lambda_m}{n^2}\right|\leq 1$. Then,
\begin{align*}
    |M_m(i\lambda_m)|&=\prod_{n=m^3}^{\infty}\left|\frac{\sin\left(\frac{i\lambda_m}{n^2}\right)}{\frac{i\lambda_m}{n^2}}\right|=\prod_{n=m^3}^{\infty}\frac{\sin(\frac{|\lambda_m|}{n^2})}{\frac{|\lambda_m|}{n^2}}\geq\prod_{n=m^3}^{\infty}\left|1-\frac{1}{6}\frac{|\lambda_m|^2}{n^4}\right|\\
    &=\exp\left(\sum_{n=m^3}^{\infty}\ln\left(1-\frac{1}{6}\frac{|\lambda_m|^2}{n^4}\right)\right)\geq\exp\left(-\frac{|\lambda_m|^2}{30}\sum_{n=m^3}^{\infty}\frac{1}{n^4}\right)\\
    &\geq\exp\left(-\frac{m^6}{30}\sum_{n=m^3}^{\infty}\frac{1}{n^2}\right)\geq\exp\left(-\frac{m^6}{30}C\right),
\end{align*}
for some $C>0$.
\end{proof}
Now we have the tools we need to construct a biorthogonal sequence to the family $(e^{\lambda_nt})_{n\in\mathbb{Z^*}}$ in
$L^2(-\frac{T}{2},\frac{T}{2})$, $T>0$.
\begin{theorem}\label{Teorema1}
There exists a constant $T_1>0$ and a biorthogonal sequence $(\Theta_m)_{m\in\mathbb{Z^*}}$ to the family $(e^{-\lambda_nt})_{n\in\mathbb{Z^*}}$ in $L^2(-\frac{T_1}{2},\frac{T_1}{2})$. Moreover,
\begin{align}\label{theta}
    \|\Theta_m\|_{L^2(-\frac{T_1}{2},\frac{T_1}{2})}\leq C\exp(bm^6),
\end{align}
where $C$ and $b$ are positive constants.
\end{theorem}
\begin{proof}
For all $m\in\mathbf{Z}^*$, let $P_m$ and $M_m$ be the functions defined in \eqref{funcion 1} and \eqref{funcion 2}, respectively. We also define the function
\begin{align*}
    \Psi_m(z)=P_m(z)\left(\frac{M_{|m|}(z)}{M_{|m|}(i\lambda_m)}\right)^{\frac{C_1}{K_1}}\frac{\sin(\delta(z-i\lambda_m))}{\delta(z-i\lambda_m)},
\end{align*}
where $\delta>0$ is an arbitrary constant, $C_1$ is given in \eqref{estPmRe} and $K_1$ in Proposition \ref{prop2}. Let
\begin{align}\label{Thetam}
    \Theta_m(t)=\frac{1}{2\pi}\int_{\mathbf{R}}\Psi_m(x)e^{itx}dx.
\end{align}
From Lemma \ref{prop1} and Proposition \ref{prop2}, we deduce that there exists $\widetilde{T}>0$, such that $\Psi_m$ is an entire function of the exponential type $\frac{\widetilde{T}}{2}$.
Moreover, from the estimates for $P_m$ and $M_m$ on the real axis (see \eqref{estPmRe} and Proposition \ref{prop2}) we obtain
\begin{equation}\label{estiTheta}
\begin{array}{l}
    \displaystyle\int_{\mathbb{R}}|\Psi_m(x)|^2dx\leq Ce^{2(2C_1+\frac{C_1K_2}{K_1})m^6}\int_{\mathbf{R}}\left|\frac{\sin(\delta(x-i\lambda_m))}{\delta(x-i\lambda_m)}\right|^2dx \\
    \qquad\qquad\qquad\,\,\leq \displaystyle\frac{C}{\delta}e^{2(2C_1+\frac{C_1K_2}{K_1})m^6}\int_{\mathbf{R}}\left|\frac{\sin t}{t}\right|^2dt\leq C_1e^{bm^6},
\end{array}
\end{equation}
where $b=2\left(2C_1+\frac{C_1K_2}{K_1}\right)$. Taking into account the properties of $\Psi_m$ and applying Paley-Wiener Theorem, we deduce that $\widehat{\Theta}_m$ has support
included in $\left(-\frac{\widetilde{T}}{2},\frac{\widetilde{T}}{2}\right)$ and $\Theta_m\in L^2(-\frac{\widetilde{T}}{2},\frac{\widetilde{T}}{2})$. Moreover, from the properties
of the inverse Fourier transform we have that the sequence $\Theta_m$ is biorthogonal to $(e^{-\lambda_mt})_{m\in\mathbf{Z}}$ in $L^2(-\widetilde{T},\widetilde{T})$. In fact,
\begin{align*}
    \int_{-\frac{\widetilde{T}}{2}}^{\frac{\widetilde{T}}{2}}\Theta_m(t)e^{\lambda_n t}dt&=\int_{-\frac{\widetilde{T}}{2}}^{\frac{\widetilde{T}}{2}}\Theta_m(t)e^{-i(i\lambda_n) t}dt=\Psi_m(i\lambda_n)
    =P_m(i\lambda_n)\frac{\sin(\delta i(\lambda_n-\lambda_m))}{\delta i(\lambda_n-\lambda_m)}=\delta_{nm}.
\end{align*}
Finally, the estimative \eqref{theta} follows from \eqref{estiTheta} by using Plancherel Theorem.
\end{proof}

\begin{remark}\label{tf} Let $\Theta_m$ be given by \eqref{Thetam}. From the proof of Theorem \ref{Teorema1}, it follows that $\widehat{\Theta}_m$ has support included
 in $(-\frac{\widetilde{T}}{2},\frac{\widetilde{T}}{2})$ and
 \begin{align}
     \|\widehat{\Theta}_m\|_{L^{\infty}(\mathbb{R})}\leq C\exp(bm^6).\nonumber
 \end{align}
\end{remark}

The following result gives the existence of a new biorthogonal sequence with better norm properties than the
one from Theorem \ref{Teorema1}. In order to prove it, for $a>0$, we define the following auxiliary functions:
\begin{align}\label{rho}
    \kappa_a=\frac{\sqrt{2\pi}}{a^2}(\chi_a*\chi_a) \quad \mbox{and}\quad \rho_m(x)=e^{x\lambda_m}\kappa_a(x),
\end{align}
where $\chi_a$ is the characteristic function of the interval
 $[-\frac{a}{2},\frac{a}{2}]$. Observe that $\kappa_a$ and $\rho_m$ satisfy the following properties:
\begin{itemize}
    \item $supp(\kappa_a)\subset[-a,a]$,
    \item $\widehat{\kappa_a}(\xi)=\frac{4}{a^2}\frac{\sin^2((\frac{a}{2})\xi)}{\xi^2}$,
    \item $\widehat{\kappa_a}(0)=1$,
    \item $supp(\rho_m)\subset[-a,a]$,
    \item $\widehat{\rho}_m(x)=\widehat{\kappa}_a(x-\lambda_m)$.
\end{itemize}

Then, we have the following result:

\begin{theorem}\label{sequebior}
There exist positive constants $T>2\pi$, $b$ and $C$ and a biorthogonal sequence $(\zeta_m)_{m\in\mathbb{Z}}$ to the family $(e^{-\lambda_mt})_{m\in\mathbb{Z}}$ in
$L^2(-\frac{T}{2},\frac{T}{2})$, with the property
\begin{align*}
    \int_{-\frac{T}{2}}^{\frac{T}{2}}\left|\sum_{n\in\mathbb{Z}^*}c_m\zeta_m(t)\right|^2dt\leq C\sum_{n\in\mathbb{Z}^*}|c_n|^2e^{2b m^6},
\end{align*}
for any sequence $(c_n)_{n\in\mathbb{N}}.$
\end{theorem}
\begin{proof}
Let $(\Theta_m)_{m\in\mathbf{Z}^*}\subset L^2(-\widetilde{T},\widetilde{T})$ be the biorthogonal sequence given by Theorem \ref{Teorema1}. Define
\begin{align*}
    \zeta_m(t)=\frac{1}{2\pi\widehat{\rho}_m(i\lambda_m)}(\Theta_m*\rho_m)(t), \ \ \ \ \ m\in\mathbb{Z}^*,
\end{align*}
where $\widehat{\rho}_m$ is the Fourier transform of $\rho_m$ defined in \eqref{rho}. Since $\zeta_m\in L^2(-\widetilde{T}-a,\widetilde{T}+a)$, take $\frac{T}{2}=\widetilde{T}+a$.
Then, applying the properties of convolution, it follows that $(\zeta_m)_{m\in\mathbb{Z}}$ is a biorthogonal sequence to $(e^{-\lambda_mt})_{m\in\mathbb{Z}}$. In fact,
\begin{align*}
    \int_{-\frac{T}{2}}^{\frac{T}{2}}\zeta_m(t)e^{\lambda_nt}dt&=\int_{-\frac{T}{2}}^{\frac{T}{2}}\zeta_m(t)e^{-i(i\lambda_n)t}dt
    =2\pi\widehat{\zeta}_m(i\lambda_n)=\frac{2\pi}{2\pi\widehat{\rho}_m(i\lambda_m)}\widehat{\Theta}_m(i\lambda_n)\widehat{\rho}_m(i\lambda_n)\\
    &=\frac{1}{\widehat{\rho}_m(i\lambda_m)}\Psi_m(i\lambda_n)\widehat{\rho}_m(i\lambda_n)=\delta_{nm}.
\end{align*}
Moreover,
\begin{align*}
    \int_{-\frac{T}{2}}^{\frac{T}{2}}\left|\sum_{m\in\mathbb{Z}^*}c_m\zeta_m(t)\right|^2dt&= \int_{-\infty}^{\infty}\left|\sum_{m\in\mathbb{Z}^*}c_m\widehat{\Theta}_m(x)\widehat{\rho}_m(x)\right|^2dx \\
    &\leq \int_{-\infty}^{\infty}\left(\sum_{m\in\mathbb{Z}^*}|c_m|\|\widehat{\Theta}_m\|_{L^{\infty}(\mathbb{R})}|\widehat{\kappa}_a(x-\lambda_m)|\right)^2dx\\
    &=\int_{-\infty}^{\infty}\left|\sum_{m\in\mathbb{Z}^*}|c_m|\|\widehat{\Theta}_m\|_{L^{\infty}(\mathbb{R})}\kappa_a(t)e^{i\lambda_mt}\right|^2dt\\
&\leq \int_{-a}^{a}\left|\sum_{m\in\mathbb{Z}^*}|c_m|\|\widehat{\Theta}_m\|_{L^{\infty}(\mathbb{R})}e^{i\lambda_mt}\right|^2dt.\\
\end{align*}
Remark that $|\lambda_{m+1}-\lambda_m|>1$, for all $m\in \mathbb{Z}^*$. Hence, from Ingham inequality and Remark \ref{tf}, we get
\begin{align}
    \int_{-a}^{a}\left|\sum_{m\in\mathbb{Z}^*}|c_m|\|\widehat{\Theta}_m\|_{L^{\infty}(\mathbb{R})}e^{i\lambda_mt}\right|^2dt\leq \sum_{m\in\mathbb{Z}^*}|c_m|^2\|\widehat{\Theta}_m\|_{L^{\infty}(\mathbb{R})}^{2}\leq\sum_{m\in\mathbb{Z}^*}|c_m|^2e^{bm^6}.
\end{align}
\end{proof}

\subsection{Controllability}

This section is devoted to prove the main result of this section. In order to do that, for any $\beta \geq b$,
where $b$ is given by Theorem \ref{sequebior},
and $f$ as in \eqref{u0,f}, we define the space
\begin{equation}\label{H1}
\widetilde{\mathcal{H}}=\left\{h\in L^2(0,2\pi):\displaystyle\sum_{k\in\mathbb{Z}}\left|\frac{\hat{h}_k}{\hat{f}_k}\right|^2e^{\beta k^6}<\infty\right\}.
\end{equation}
Then, our main result reads as follows:
\begin{theorem}\label{controla}
Let $f\in L^2(0,2\pi)$ a function verifying \eqref{u0,f} and $\widetilde{\mathcal{H}}$ defined by \eqref{H1}. There exists $T>0$, such that, for any initial data $u_0\in \widetilde{\mathcal{H}}$,
there exist a control $v\in L^2(0,T)$ for which the solution of \eqref{a01} satisfies $u(T,x)=0$.
\end{theorem}

\begin{proof}
Let $T>2\pi$ and $(\zeta_m)_{m\in\mathbb{Z}^*}$  given by Theorem \ref{sequebior}. For $u_0\in \widetilde{\mathcal{H}}$,
such that $u_0(x)=\displaystyle\sum_{k\in \mathbf{Z}^*}\widehat{u}_{k}^{0}e^{ikx}$, define $v$ as follows:
\begin{align}
    v(t)=-\sum_{m\in\mathbb{Z}^*}\frac{\hat{u}_{m}^{0}}{\hat{f}_m}e^{\frac{T}{2}\lambda_m}\zeta_m\left(t-\frac{T}{2}\right),\ \ \ \ t\in(0,T).
\end{align}
From the properties of the biorthogonal sequence $(\zeta_m)_{m\in\mathbb{Z}}$, we deduce that $v$ is a control that satisfies \eqref{equiv}, i. e., leads the solution to zero.
Moreover, $v\in L^2(0,T)$. In fact, from Theorem \ref{sequebior},
\begin{align}
    \int_{0}^{T}|v(t)|^2dt&=
    \int_{0}^{T}\left|-\sum_{m\in\mathbb{Z}^*}\frac{\hat{u}_{m}^{0}}{\hat{f}_m}e^{\frac{T}{2}\lambda_m}\zeta_m\left(t-\frac{T}{2}\right)\right|^2dt
    \leq C\sum_{m\in\mathbb{Z}^*}\frac{|\hat{u}_{m}^{0}|^2}{|\hat{f}_m|^2}e^{bm^6}\leq C ||u_0||_{\widetilde{\mathcal{H}}}^2,
\end{align}
for some $C>0$.
\end{proof}

\section{The Nonlinear System}

This section is devoted do analyze the controllability of the full system \eqref{a01-1}. Our main result reads as follows:

\begin{theorem}\label{controla-nl}
Let $f\in L^2(0,2\pi)$ a function verifying \eqref{u0,f} and $\widetilde{\mathcal{H}}$ defined by \eqref{H1}. There exists $T>0$ and $\delta>0$,
such that, for any $u_0,u_1\in \mathcal{H}$ satisfying
$$\|u_0\|_{\mathcal{H}}\leq\delta\ ,\ \|u_1\|_{\mathcal{H}}\leq\delta,$$
there exist a control $v\in L^2(0,T)$, such that system \eqref{a01-1} admits a solution $u\in C([0,T];\widetilde{\mathcal{H}})$ verifying
$$u(0,x)=u_0\,\ \ u(T,x)=u_1.$$
\end{theorem}

In order to prove Theorem \ref{controla-nl}, we combine the analysis of the linearized system, a fixed point argument and some Bourgain smoothing properties of the Kawahara equation on a periodic domais. Therefore, some technical results are needed.

We start by introducing the Bourgain spaces associated to  the Kawahara equation on $\mathbb{T}$. We remark that this is equivalent to impose the periodic boundary conditions over the interval $(0, 2\pi)$, as in \eqref{a01-1}.

For given $b,s\in \mathbb{R}$ and a function $u : \mathbb{R}\times\mathbb{T} \rightarrow \mathbb{R}$, we define the spaces
\begin{align*}
\|u\|_{X_{b,s}^{\beta,\gamma}}=\left(\sum_{k\in\mathbb{Z}}\int_{\mathbb{R}}\langle k\rangle^{2s}\langle \tau-p^{\beta,\gamma}(k)\rangle^{2b}|\widehat{u}(\tau,k)|^2d\tau\right)^{\frac{1}{2}}
\end{align*}
and
\begin{align*}
\|u\|_{Y_{b,s}^{\beta,\gamma}}=\left(\sum_{k\in\mathbb{Z}}\left(\int_{\mathbb{R}}\langle k\rangle^{s}\langle \tau-p^{\beta,\gamma}(k)\rangle^{b}|\widehat{u}(\tau,k)|d\tau\right)^2\right)^{\frac{1}{2}},
\end{align*}
where $\widehat{u}(\tau,k)$ denotes the Fourier transform of $u$ with respect to the time variable $t$ and the space variable $x$, $\langle \cdot\rangle=\sqrt{1+|\cdot|^2}$ and $p^{\beta,\gamma}(k)=\beta k^3-\gamma k.$

The spaces $X_{b,s}$ and $Y_{b,s}$ are the completion of the Schwartz space $\mathcal{S}(\mathbb{R}\times\mathbb{T})$ under the norm $\|u\|_{X_{b,s}}$ and $\|u\|_{Y_{b,s}}$, respectively. Observe that, for any $u\in X_{b,s}$,
\begin{align*}
    \|u\|_{X_{b,s}}=\|S(-t)u\|_{H^b(\mathbb{R},H^s(\mathbb{T}))}.
\end{align*}
For given $b,s\in\mathbb{R}$, let us introduce the space
\begin{align*}
    Z_{b,s}=X_{b,s}\cup Y_{b-\frac{1}{2},s}
\end{align*}
endowed with the norm
\begin{align*}
 \|u\|_{Z_{b,s}}= \|u\|_{X_{b,s}}+ \|u\|_{Y_{b-\frac{1}{2},s}}.
\end{align*}
For a given interval $I$, we denote by $X_{b,s}(I)$ and $Z_{b,s}(I)$ the restriction spaces $X_{b,s}$ to the interval $I$ with the norms
\begin{align*}
 \|u\|_{X_{b,s}(I)}=\inf\{\|\widetilde{u}\|_{X_{b,s}}|\widetilde{u}=u\ \ \text{on}\ \ \mathbb{T}\times I\}\,\mbox{ and }\,
 \|u\|_{Z_{b,s}(I)}=\inf\{\|\widetilde{u}\|_{Z_{b,s}}|\widetilde{u}=u\ \ \text{on}\ \ \mathbb{T}\times I\}.
\end{align*}
If $I=(0,T)$, for simplicity, we denote $X_{b,s}(I)$ and $Z_{b,s}(I)$ by $X_{b,s}^T$ and $Z_{b,s}^T$, respectively.
Finally, let
\begin{align*}
    \mathbb{Z}_{\frac{1}{2},s}^{T}=Z_{\frac{1}{2},s}^{T}\cup C([0,T];H^s(\mathbb{T})).
\end{align*}

The following estimates related to the Bourgain space $X_{b,s}^T$ and $Z_{b,s}^T$ play
important roles in the proof of Theorem \ref{controla-nl}.
\begin{lemma}\label{bour1}
Let $b,s\in\mathbb{R}$ and $T>0$ given. There exist a constant $C>0$ such that
\begin{enumerate}
    \item For any $\varphi\in H^s(\mathbb{T})$,
\end{enumerate}
\begin{align*}
    \|S(t)\varphi\|_{\mathbb{Z}_{\frac{1}{2},s}^T}\leq C\|\varphi\|_{s}
\end{align*}
\item For any $f\in Z_{-\frac{1}{2},s}^T$,
\begin{align*}
    \left\|\int_{0}^{t}S(t-\tau) f(\tau)d\tau\right\|_{\mathbb{Z}_{\frac{1}{2},s}^T}\leq C \|f\|_{Z_{-\frac{1}{2},s}^T}.
\end{align*}
\end{lemma}
\begin{proof}
See \cite{hirayama}.
\end{proof}
\begin{lemma}\label{biliest}
    Let $s\geq -1$ and $T>0$ be given. There exist a constant $C$ such that the following bilinear estimate
    \begin{align*}
        \|(uv)_x\|_{Z_{-\frac{1}{2},s}^{T}}\leq C \|u\|_{Z_{\frac{1}{2},s}^{T}}\|uv\|_{Z_{\frac{1}{2},s}^{T}}
    \end{align*}
    holds.
\end{lemma}
\begin{proof} See \cite{hirayama} and Lemma 3.2 in \cite{zhang1}.
\end{proof}
From now on, we can pass to the proof of Theorem \ref{controla-nl}.
\begin{proof} Throughout the proof we consider the following space defined above:
\begin{align*}
    \mathbb{Z}_{\frac{1}{2},s}^{T}=Z_{\frac{1}{2},s}^{T}\cup C([0,T];\widetilde{\mathcal{H}})).
    \end{align*}

    Given $u_0,u_1\in\widetilde{\mathcal{H}}$ and $\xi\in C([0,T;\widetilde{\mathcal{H}}])$, let
    $$\zeta^0=S(T)u_0\ \ \text{and}\ \ \zeta^1=\int_{0}^{T}S(T-\tau)(\xi\xi_x)(\tau)d\tau.$$
    According to Theorem \ref{controla}, there exists $v\in L^2(0,T),$ such that
    \begin{align}\label{I}
        \varphi(t)=\int_{0}^{t}S(t-\tau)f(x)v(\tau)d\tau
    \end{align}
    satisfies
   \begin{align}\label{II}
       \varphi(0)=0\ \ \text{and}\ \ \varphi(T)=u_1-\zeta^0-\zeta^1.
   \end{align}
   Moreover,
   $$\|v\|_{L^2(0,T)}\leq C(\|u_0\|_{\widetilde{\mathcal{H}}}+\|u_1\|_{\widetilde{\mathcal{H}}}+\|\zeta^1\|_{\widetilde{\mathcal{H}}}).$$
   Then, we can define a nonlinear map $\Psi:C([0,T];\widetilde{\mathcal{H}})\rightarrow L^2(0,T)$ as follows:
   $$\Psi(\xi):=h,$$
   where $h$ is a control verifying \eqref{I} and \eqref{II}.

   With the notation introduced above, we define the following nonlinear map $\Gamma$ from $\mathbb{Z}_{\frac{1}{2},0}^T$ into itself:
    \begin{align}\label{duha}
        \Gamma (\xi)(t)=S(t)u_0+\int_{0}^{t}S(t-\tau)[f(x)\Psi(\xi))(\tau)+(\xi\xi_x)(\tau)]d\tau.
    \end{align}
If we can prove that $\Gamma$ is a contraction map, then its fixed point $u$ is a solution of \eqref{a01-1}, with $h=\Psi(u)$ and
satisfies $u(T,x)=u_1$.

    Applying Lemmas \ref{bour1} and \ref{biliest}, we obtain the following estimate:
    \begin{align*}
        \|\Gamma (\xi)\|_{\mathbb{Z}_{\frac{1}{2},0}^T}&\leq C\|u_0\|_{\widetilde{\mathcal{H}}}+C\left\|\int_{0}^{t}S(t-\tau)[f(x)\Psi(\xi))(\tau)+(\xi\xi_x)(\tau)\right\|_{\mathbb{Z}_{\frac{1}{2},0}^T}\\
        &+C\| \xi\xi_x\|_{\mathbb{Z}_{-\frac{1}{2},0}^T}\\
        &\leq C\|u_0\|_{\widetilde{\mathcal{H}}}+C\left\|f(x)\Psi(\xi)(\tau)\right\|_{L^2(0,T; L^2(\mathbb{T}))}+C\|\xi\|_{\mathbb{Z}_{\frac{1}{2},0}^T}^{2}\\
        &\leq C\|u_0\|_{\widetilde{\mathcal{H}}}+C(\|u_1\|_{\widetilde{\mathcal{H}}}+\|u_0\|_{\widetilde{\mathcal{H}}}+\|\zeta^1\|_{\widetilde{\mathcal{H}}})
        +C\|\xi\|_{\mathbb{Z}_{\frac{1}{2},0}^T}^{2}.
    \end{align*}
 Observe that
    \begin{align*}
        \|\zeta^1\|_{\widetilde{\mathcal{H}}}=\left\|\int_{0}^{T}S(T-\tau)(\xi\xi_x)(\tau)d\tau\right\|_{\widetilde{\mathcal{H}}}
        \leq C \sup_{t\in[0,T]}\left\|\int_{0}^{t}S(t-\tau)(\xi\xi_x)(\tau)d\tau\right\|_{\widetilde{\mathcal{H}}}
        \leq C\|\xi\|_{\mathbb{Z}_{\frac{1}{2},0}^T}^{2}.
    \end{align*}
    Hence,
    \begin{align*}
         \|\Gamma (\xi)\|_{\mathbb{Z}_{\frac{1}{2},0}^T}\leq C(\|u_1\|_{\widetilde{\mathcal{H}}}+\|u_0\|_{\widetilde{\mathcal{H}}})+C\|\xi\|_{\mathbb{Z}_{\frac{1}{2},0}^T}^{2}.
    \end{align*}
    For $R>0$, let $B_R$ be a bounded subset of $\mathbb{Z}^{T}_{\frac{1}{2},0}$:
    $$B_R=\{g\in\mathbb{Z}^{T}_{\frac{1}{2},0}\ | \ \|g\|_{\mathbb{Z}_{\frac{1}{2},0}^T}\leq R\}.$$
    We choose $\delta>0$ and $R>0$, such that
    $$2C\delta+CR^2\leq R,\ \ \ CR\leq \frac{1}{2}.$$
    Then, $\|\Gamma (u)\|_{\mathbb{Z}_{\frac{1}{2},0}^T}\leq R$, that is, $\Gamma$ map $B_R$ into itself. In addition, for any $u, v\in B_R$,
similarly, we have
$$\|\Gamma(u)-\Gamma(v)\|_{\mathbb{Z}_{\frac{1}{2},0}^T}\leq\frac{1}{2}\|u-v\|_{\mathbb{Z}_{\frac{1}{2},0}^T}.$$
$\Gamma$ is thus a contracting map on $B_R$. By the Banach fixed point theorem, there is a
unique solution to the integral equation \eqref{duha} which is the desired solution of \eqref{a01-1}.
\end{proof}

\section{Comments and Open Problems}

We close this paper with some comments and open problems that are worthy of further study:

\begin{itemize}

\item In \cite{micu-ortega-pazoto}, the authors consider the following parabolic type control system
\begin{equation}
\begin{cases}
u_t + i(-\partial_{xx}^2)^{\frac{1}{2}}u-\varepsilon\partial_{xx}^2 u=f(x)v_\varepsilon (t), &\text{ in  } \ \ (0,T)\times (0,\pi),\\
u(t,0)=u(t,\pi)=0 &\text{ in  } \ \ (0,T),\\
u(0,x)=u_0(x), &\text{ in  } \ \ (0,\pi),\nonumber
\end{cases}
\end{equation}
where $v_\varepsilon$ is a control and $f$ is a given profile. For $\varepsilon=0$ the system is of hyperbolic type and the authors
show that the control steering the hyperbolic system to rest  can be approximated by a sequence $(v_\varepsilon)_{\varepsilon > 0}$
of controls of the parabolic system when $\varepsilon\rightarrow 0$. The proof is based on the moment problem with respect to the
nonharmonic Fourier family $(e^{\lambda_n})_{n\in \mathbb{N}}$, where $\lambda_n=in-\varepsilon n^2$, $n\geq 1$, are the eigenvalues of the
corresponding differential state operator. More recently, in \cite{micu-bugariu}, the same problem was studied for the linear wave equation
by introducing a viscous term which contains a fractional power of the Dirichlet Laplace operator.
It is a difficult problem that remains unanswered for the Kawahara equation.

\item  Employing the same approach, Theorems \ref{controla} and Theorems \ref{controla-nl} can be proved for
the KdV equation with similar statements.
In this case, our analysis can be simplified due to the absence of the fifth order dispersive term.

\item Taking into account the results obtained in \cite{yan}, we expect that our analysis can be extended for the modified
Kawahara equation. Moreover, other types of controls could be considered, such as boundary or moving controls. We refer to the works \cite{rosier,rosier-zhang2,russell-zhang2} in which the control problem was addressed in the context of the Korteweg-de Vries and the Benjamin-Bona-Mahony equations.

\end{itemize}

\section*{Acknowledgements} The first author was partially support by CNPq (Brazil). The second author was suppported by CAPES and CNPq (Brazil).

\section*{Declarations}

\noindent {\bf Ethical Approval:} Not applicable.

\noindent {\bf Competing interests:} The authors have no competing interest to declare.

\noindent {\bf Authors' contributions:} All authors contributed equally to the final manuscript.

\noindent {\bf Funding:} Not applicable.

\noindent {\bf Data availability:} There is no data associated with this manuscript.

\end{document}